\documentstyle{article}     

\newtheorem{th}{Theorem}
\newtheorem{ax}{Axiom}
\newtheorem{lm}{Lemma}
\newtheorem{df}{Definition}
\newtheorem{pr}{Proposition}
\newtheorem{cl}{Corollary}
\newtheorem{as}{Assumption}
\newcommand{\bth}{\begin{th}\hspace{-5pt}{\bf .} \ }
\newcommand{\eth}{\end{th}}
\newcommand{\bax}{\begin{ax}\hspace{-5pt}{\bf .} \ }
\newcommand{\eax}{\end{ax}}
\newcommand{\blm}{\begin{lm}\hspace{-5pt}{\bf .} \ }
\newcommand{\elm}{\end{lm}}
\newcommand{\bdf}{\begin{df}\hspace{-5pt}{\bf .} \ }
\newcommand{\edf}{\end{df}}
\newcommand{\bpr}{\begin{pr}\hspace{-5pt}{\bf .} \ }
\newcommand{\epr}{\end{pr}}
\newcommand{\bcl}{\begin{cl}\hspace{-5pt}{\bf .} \ }
\newcommand{\ecl}{\end{cl}}
\newcommand{\bas}{\begin{as}\hspace{-5pt}{\bf .} \ }
\newcommand{\eas}{\end{as}}

\newcommand{\bit}{\begin{itemize}}
\newcommand{\eit}{\end{itemize}\par\noindent}
\newcommand{\beq}{\begin{equation}}
\newcommand{\eeq}{\end{equation}\par\noindent}
\newcommand{\beqa}{\begin{eqnarray*}}
\newcommand{\eeqa}{\end{eqnarray*}\par\noindent}
\newcommand{\beqn}{\begin{eqnarray}}
\newcommand{\eeqn}{\end{eqnarray}\par\noindent}

\newcommand{\ela}{\end{array}\right.}
\newcommand{\bra}{\left\{\begin{array}{r}}
\newcommand{\era}{\end{array}\right.}

\newcommand{\ben}{\begin{enumerate}}
\newcommand{\een}{\end{enumerate}\par\noindent}

\begin{document}
\
\par\medskip\par\noindent  
\centerline{\large{\bf DISJUNCTIVE QUANTUM LOGIC}}
\par\smallskip\par\noindent
\centerline{\large{\bf IN DYNAMIC PERSPECTIVE}}
\par\medskip\par\noindent
\centerline{\normalsize{B{\scriptsize OB} C{\scriptsize OECKE}}}
\par\medskip\par\noindent
\centerline{\small{Free University of Brussels, Department of Mathematics,}}
\par\noindent
\centerline{\small{Pleinlaan 2, B-1050 Brussels\,; bocoecke@vub.ac.be\,;}}
\par\smallskip\par\noindent
\centerline{\small{\sl and}}
\par\smallskip\par\noindent
\centerline{\small{Imperial College of Science, Technology \&
Medicine, Theoretical Physics Group,}}
\par\noindent 
\centerline{\small{The Blackett Laboratory, South Kensington, London SW7 2BZ\,.}}
\par\bigskip\par\noindent 
\centerline{\small{{\it Current address:}}}

\smallskip\noindent
\centerline{\small{University of Oxford, Computing 
Laboratory\,,  Wolfson Building, Parks Road,}} 

\noindent 
\centerline{\small{Oxford, OX1 3QD, UK\,; e-mail:
coecke@comlab.ox.ac.uk\,.}} 
\begin{abstract}
\noindent
In Coecke (2002a) we proposed the {\sl intuitionistic or disjunctive
representation of quantum
logic}, i.e., a representation of the property lattice of physical
systems as a complete
Heyting algebra of logical propositions on these properties, where
this complete Heyting algebra goes equipped with an
additional operation, the {\sl operational resolution}, which
identifies the properties within the logic of
propositions.  This representation has an important application
``towards dynamic quantum logic'', namely in describing the
temporal indeterministic propagation of actual properties of physical
systems.  This paper can as such by conceived as
an  addendum to ``Quantum Logic in Intuitionistic Perspective'' that
discusses spin-off and thus provides an additional
motivation.  We derive a quantaloidal semantics for dynamic disjunctive quantum
logic and illustrate it for the particular case of a perfect (quantum) measurement. 
\end{abstract}
\par
\medskip
\par
\noindent {Key words: Quantum logic, dynamic logic, intuitionistic
logic, property lattice, operational resolution, quantaloid.}

\bigskip\noindent
{\bf 1. INTRODUCTION}

\medskip\noindent
In Amira, Coecke and Stubbe (1998), Coecke and Stubbe (1999), Coecke
(2000), Coecke, Moore and Stubbe (2001),
Coecke and Smets (2000) and Sourbron (2000) steps have been taken
towards a dynamic quantum logic, to great extend inspired by
the representation theorem for Schr\"odinger flows of Faure, Moore
and Piron (1995), which itself incorporates the
results of Faure and Fr\"olicher  (1993, 1994) on categorical representations of
projective geometries; for previous attempts in that direction we refer to Pool (1968)
and Daniel (1989).  The crucial formal notion in this new approach is
that of an {\sl operational resolution}, a map that
assigns to collections of either states or properties of a physical
system the strongest property whose actuality is implied by
that of each member  in the collection\,: Indeterministic state
transitions and property transitions are then exactly described
by those maps between powersets of either the state space or the
property lattice that preserve the operational resolution.
Formally, this discussion takes place in the category of so called
{\sl quantaloids} and {\sl quantaloid morphisms}, i.e., the
category of {\it sup}-{\sl lattice enriched categories}. In this paper an ad hoc
definition is given.

\smallskip
The notion of an operational resolution also emerged in a different
context\,: If one represents the
lattice of properties of a physical system in terms of {\sl logical
propositions on these properties}, the operational
resolution comes in as an additional operation which identifies
physical properties within this propositional logic, and which
moreover establishes this representation as a true equivalence
(Coecke 2002a, Section 4). However, the domain of the operational
resolution in Coecke (2002a) is a restriction of the one introduced
in Coecke and Stubbe (1999)\,: it is not the powerset of
the property lattice  but its Bruns-Lakser distributive hull (Bruns
and Lakser 1970), this since it is the latter that
constitutes the logical propositions with respect to actuality of
physical properties. The main message of this paper will as
such be  ``imposing a refinement on the operational resolution in the
capacity as the mathematical object that generates
state and property transitions, guided by logical analysis of its domain''.

\smallskip
We refer to Coecke (2002a) for the preliminaries to this paper on states,
properties,  actuality of properties, Cartan maps and actuality
sets\,; superposition states, superposition properties and
superpositional faithfulness\,;  atomistic lattices, complete
lattices, Galois adjoints, complete Heyting algebras,
Bruns-Lakser distributive ideals, Bruns-Lakser distributive hulls and
distributive join dense closures; ortholattices,
orthomodular lattices and Sasaki projections; for the latter we also
refer to Piron (1976) and Kalmbach (1983)\,. For a brief
outline of ordinary and enriched category theory we refer to Borceux
and Stubbe (2000), for a detailed one to Borceux (1994).
The particular case of quantaloids is discussed in
Rosenthal (1991).

\bigskip\noindent
{\bf 2. TOWARDS A DYNAMIC QUANTUM LOGIC}

\medskip\noindent
In Amira, Coecke and Stubbe (1998) it was shown that the inducible
state and property transitions on a physical
system, the procedures that realize these transitions being called
{\sl inductions}, constitute a quantale, i.e.,
a complete lattice $\bigl(L,\bigvee\bigr)$ equipped with an
additional associative operation $-\&-:L\times L\to
L$ that distributes over suprema at both sides\,:
$$
a\&\Bigl(\bigvee B\Bigr)=\bigvee_{b\in B}(a\& b)
\quad\quad\quad\quad
\Bigl(\bigvee A\Bigr)\&b=\bigvee_{a\in A}(a\& b)\,.
$$
Note here that the collection of all inductions that can be
effectuated on a particular physical system include both
measurements and evolution, the latter to be understood as ``let the
system evolve''. To fix ideas, let us consider the example
of a {\sl perfect measurement induction}
$e_{PM}$ as it is outlined in Coecke and Smets (2000), based on the
notion of perfect measurement in Piron
(1976)\,.\footnote{These perfect measurements encode in quantum
logical terms the
measurements in orthodox quantum theory represented by self-adjoint
operators with a binary spectrum.} Given a system described
by a complete orthomodular lattice $L$\,, e.g., a classical system or
a quantum system, then actuality of a property $a\in L$ in
a measurement with as {\sl eigenproperties} $b$ and $b'$  guarantees
actuality of {\it either}
$$
\varphi_b(a):=b\wedge(a\vee b')
\quad\quad{\it or}\quad\quad
\varphi_{b'}(a):=b'\wedge(a\vee b)\,,
$$
i.e., the Sasaki projection of $a$ on either $b$ or 
$b'$\,.\footnote{Recall here, as mentioned in
Footnote 8 of Coecke
(2002a), that we consider a quantum measurement as an external action
on the system that changes its state, and as such also
its actual properties.} Note that in the case that one of the two alternatives turns
out to be $0$ the outcome is determined (since $0$ is impossible). If this lattice is
moreover atomistic and satisfies the covering law, what is still the case both for
classical and quantum systems (Piron 1976), then, provided that the states are
encoded as atoms,
$e_{PM}$ imposes a change of the initial state to {\it either}
$$
\varphi_b(p):=b\wedge(p\vee b')
\quad\quad{\it or}\quad\quad
\varphi_{b'}(p):=b'\wedge(p\vee b)\,,
$$
again provided that none of the alternatives yields $0$\,, in which case we only
consider the non-$0$ outcome. To
$e_{PM}$ we can as such attribute a map
$$
\tilde{e}_{PM}:{\cal P}(\Sigma)\to{\cal P}(\Sigma):
T\mapsto\left\{b\wedge(p\vee b'),b'\wedge(p\vee
b)\bigm|p\in T\right\}\setminus\{0\}
$$
that assigns the possible outcome states whenever the system is
initially in a state in $T$\,.
For the more general case of a property transition the situation is however
somewhat more complicated to describe.  Although at first sight one
could go for a union
preserving map between powersets of the property lattice as it is
done in Coecke and Smets (2000)\,,
saying that actuality of $a$ guarantees actuality of
{\it either} $\varphi_b(a)$ {\it or} $\varphi_{b'}(a)$ is indeed
somewhat ambiguous. Whenever
$\varphi_b(a)$ is actual also any $c\geq\varphi_b(a)$ is actual, so
one might additionally want to elucidate something in the
sense of ``one focuses on maximally strong possible outcomes'',
whatever this might mean. But then again, taking for example an
atomistic property lattice where states are encoded as atoms,
actuality of $\varphi_b(a)$ also guarantees that some state in
${\{p\in\Sigma|p\leq\varphi_b(a)\}}$ is actual, what illustrates that
the above elucidation is indeed sloppy.
The key to solve this problem is exactly the {\sl logic of actuality
sets} proposed in Coecke (2002a)\,, a presentation of the
logical propositions on properties of a physical system in terms of
actuality, which indicates that, under the assumption of
superpositional faithfulness we have to define
$$
\hat{e}_{PM}:{\cal DI}(L)\to{\cal DI}(L): A\mapsto{\cal
C}\Bigl(\left\{b\wedge(a\vee
b'),b'\wedge(a\vee b)\bigm|a\in A\right\}\Bigr)
$$
as the map that describes propagation of actuality sets in a perfect
measurement,
recalling here that ${\cal C}$ stands for the composite of the
implicative and disjunctive
closure, i.e., 
$$
{\cal C}(A):=\bigl\{\bigvee_LB\bigm|B\subseteq\downarrow\![A]\cap{\cal
D}(L)\bigr\}\,.
$$ 
Turning back the
collection of all inductions that can be effectuated on a particular
physical system, this collection
being denoted as ${\cal E}$\,, we obtain as such two quantales
$$
\tilde{\cal E}:=\left\{\tilde{e}:{\cal P}(\Sigma)\to{\cal
P}(\Sigma)\bigm|e\in{\cal E}\right\}
$$
$$
\hat{\cal E}_{\cal DI}:=\left\{\hat{e}:{\cal DI}(L)\to{\cal
DI}(L)\bigm|e\in{\cal E}\right\}
$$
respectively expressing the inducible state and property transitions
for this system, provided
the system is not destroyed.\footnote{It is merely a technicality to
avoid the restriction ``provided
the system is not destroyed''\,; for details on this matter we refer
to Coecke, Moore and Stubbe (2001).}
The suprema in these quantales are  calculated pointwisely with respect to the
suprema of the codomain, the quantale multiplication coincides with
composition of maps, and closure of these quantales under
these operations is to be understood in terms of inductions
respectively as ``arbitrary choice on effectuation'' and
``consecutive effectuation'' (Amira, Coecke and Stubbe 1998\,, Coecke
and Stubbe 1999)\,. Once at this point, one might prefer
to express
$\hat{\cal E}_{\cal DI}$ rather in terms of the distributive hull $H$
of $L$ than in terms of actuality sets, i.e.,
in terms of {\sl propagation of logical propositions on
properties} rather than in terms of {\sl propagation of
actuality sets}\,:
$$
\hat{\cal E}:=\left\{\hat{e}:H\to H\bigm|e\in{\cal E}\right\}\,,
$$
recalling from Bruns and Lakser (1970) that the inclusion
$i:L\hookrightarrow H$ satisfies
\[\begin{array}{ccc}
\ \,L&\stackrel{i}{\hookrightarrow}&H\ \,\\
{\scriptstyle\cong}\updownarrow\ & &\ \updownarrow{\scriptstyle\cong}\\
\downarrow\![L]&\hookrightarrow&{\cal DI}(L)
\end{array}\]
where the isomorphism between ${\cal DI}(L)$ and $H$ is
established via $A\mapsto\bigvee_H A$ and that between
$\downarrow\![L]$ and $L$ via $A\mapsto\bigvee_L A$\,.

\smallskip 
It makes however sense to generalize the above to maps where the
codomain is different from the domain, and this for two
reasons\,: (i) It allows to describe ``change of system'', where we
conceive a system as being
defined exactly by its set of distinct (with respect to the
corresponding actual properties) possible
realizations, i.e., by its set of states\,; (ii) Besides temporal
propagation it also allows to encode entanglement, or any
form of interaction including separation, in terms of ``mutual
induction of properties'' (Coecke 2000)\,.
Therefore we will extend our framework from quantales to quantaloids,
i.e., categories enriched in {\it sup}-lattices. Explicitly, a {\sl
quantaloid} $\underline{{\rm Q}}$ is a category in which all the
morphism sets $\underline{{\rm Q}}(A,B)$ are complete lattices with the
ordering such that the by $f:A\to B$ induced morphism actions
$f\circ-:\underline{{\rm Q}}(B,C)\to\underline{{\rm Q}}(A,C)$
and $-\circ f:\underline{{\rm Q}}(C,A)\to\underline{{\rm Q}}(C,B)$ preserve
suprema. Quantaloid morphisms are those
functors that preserve suprema when restricted to morphism sets
(Rosenthal 1991).
Quantales with multiplicative unit are then exactly one-object quantaloids. The unit
in our setting is provided by the induction ``freeze'' with obvious significance.
Denoting the quantaloid of complete lattices and
{\it sup}-morphisms as $\underline{{\rm Sup}}$ we then have that
$\tilde{\cal E}\hookrightarrow\underline{{\rm Sup}}\bigl({\cal
P}(\Sigma_1),{\cal P}(\Sigma_2)\bigr)$ and $\hat{\cal E}_{\cal
DI}\cong\hat{\cal
E}\hookrightarrow\underline{{\rm Sup}}(H_1,H_2)\cong\underline{{\rm Sup}}\bigl({\cal
DI}(L_1),{\cal DI}(L_2)\bigr)$ are functorial {\it
sup}-inclusions, where functorial is to be understood in the sense
that any composition of inductions encodes as
$\underline{{\rm Sup}}$-composition. However, as we will see below, there is
an additional feature to this inclusion.

\vfill\eject\noindent
{\bf 3. CAUSAL CONTINUITY}

\medskip\noindent
Clearly, there is a strong analogy of the above with the {\sl
non-commutative geometric logic} or {\sl observational
semantics} of Abramsky and Vickers (1993) and Resende (2000) that has
been developed in order to describe sequences of
interaction with and observation of computational devices.  However,
as for example mentioned in Resende (2000
\S 3.1), the observational semantics proposed in Abramsky and Vickers
(1993) is not
applicable to quantum processes, this in particular since in quantum
processes both the suprema in the property lattice
and disjunctions of properties are essential. We will now show how
this implements formally
within the above setting.

\smallskip
We will refer by {\sl propagation of strongest actual properties}
with respect to an induction $e$ to the map $\bar{e}:L_1\to
L_2$ that assigns to a property $a\in L_1$ the strongest property
$b\in L_2$ of which actuality after effectuating $e$ is
guaranteed by actuality of $a$ before effectuating $e$\,.
Following  Faure, Moore and Piron (1995), Coecke (2000) and Coecke,
Moore and Stubbe (2001),
the map that describes propagation of strongest actual properties
preserves suprema.
This follows from the fact that propagation is adjoint to causal
assignment (Coecke 2000\,, Coecke, Moore and Stubbe 2001).  Roughly,
the argument goes as follows\,: By
conjunctivity of infima in property lattices it follows that {\sl
assignment of weakest causes of actuality} preserves all
non-empty infima,  and as such, it induces a unique left Galois
adjoint on the upper pointed extensions of the involved property lattices
(Coecke and Moore 2000\,, Coecke, Moore and Stubbe
2001)\,;\footnote{The formal need to consider these upper pointed
extensions formally implements the physical possiblity of destruction of the
system, and as such embodies a way to
avoid the restriction ``provided the system is not destroyed''
mentioned in the previous footnote.  A discussion concerning can also be found in
Sourbron (2000).} one then verifies that this left Galois adjoint expresses
propagation of strongest actual properties.

\smallskip 
Since, given an actuality set $A$\,, the strongest property
that is actual with certainty is exactly $\bigvee A$\,, a
role that is played for states by the operational resolution
${\cal R}_\Sigma:{\cal P}(\Sigma)\to L$ sensu Coecke and
Stubbe (1999) and discussed in the first paragraph of the
introduction to this paper, it then follows that for any map
in
$\tilde{\cal E}$ or $\hat{\cal E}_{{\cal DI}}$ there  exists
a map
$\bar{e}:L_1\to L_2$ such that we respectively have
commutation of
\[
\begin{array}{ccc}
L_1&\stackrel{\bar{e}}{\longrightarrow}&L_2\\
\hspace{-3mm}{\scriptstyle{\cal R}_{\Sigma_1}}\uparrow\ \ &
&\
\uparrow{\scriptstyle{\cal R}_{\Sigma_2}}\hspace{-5mm}\\
{\cal
P}(\Sigma_1)&\stackrel{\tilde{e}}{\longrightarrow}&{\cal
P}(\Sigma_2)
\end{array}
\quad\quad\quad\quad
\begin{array}{ccc}
\ \ L_1&\stackrel{\bar{e}}{\longrightarrow}&L_2\ \ \\
\hspace{-3mm}{\scriptstyle\bigvee_1(-)}\uparrow\ \ & &\
\uparrow{\scriptstyle\bigvee_2(-)}\hspace{-5mm}\\ {\cal
DI}(L_1)&\stackrel{\hat{e}_{{\cal
DI}}}{\longrightarrow}&{\cal DI}(L_2)
\end{array}
\] which translates in terms of $\hat{\cal E}$ as commutation
of
\[
\begin{array}{ccc}
L_1&\stackrel{\bar{e}}{\longrightarrow}&L_2\,\\
\hspace{-3mm}{\scriptstyle{\cal R}_1}\uparrow\ \ &
&\uparrow{\scriptstyle{\cal R}_2}\ \hspace{-5mm}\\
H_1&\stackrel{\hat{e}}{\longrightarrow}&H_2
\end{array}  
\] where we slightly abused notation by restricting the
codomain of the operational resolution. Note that when
replacing in the above ${\cal DI}(L)$ by ${\cal P}(L)$, i.e.,
requiring for a union preserving map $g:{\cal P}(L_1)\to{\cal
P}(L_2)$ that there exists a map $f:L_1\to L_2$ such that
$\bigvee_2\bigl(g(-)\bigr)=f\bigl(\bigvee_1(-)\bigr)$ does
not assure existence of a map
$h:{\cal DI}(L_1)\to{\cal DI}(L_2)$ such that we have
commutation of
\beq\label{eq:comm}
\begin{array}{ccc} {\cal
DI}(L_1)&\stackrel{h}{\longrightarrow}&{\cal DI}(L_2)\\
\hspace{-3mm}{\scriptstyle{\cal C}_1}\uparrow\ \ & &\
\uparrow{\scriptstyle{\cal C}_2}\hspace{-5mm}\\ {\cal
P}(L_1)&\stackrel{g}{\longrightarrow}&{\cal P}(L_2)
\end{array}\eeq It suffices to note that ${\cal
DI}(L)$-suprema and ${\cal P}(L)$-suprema don't coincide.
Therefore, the considerations made in this paper reveal this
aspect as an additional feature of the maps in $\hat{\cal E}$
on propagation of actuality sets besides the one imposed by
preservation of suprema for propagation of strongest actual
properties.  In particular can all this be encoded as the
factorization of quantaloid morphisms expressed in the
following commutative diagram in the category of
quantaloids\,:
\[
\begin{array}{ccccc} &&&\hspace{-6mm}\underline{{\rm
Sup}}\hspace{-5mm}\\ &&^G \nearrow& &\nwarrow^H\\
\hat{\cal E}&\hookrightarrow&\underline{{\rm
DCHeyt}}&\stackrel{F}{\leftarrow}&\underline{{\rm PSup}}
\end{array}
\] where
\begin{itemize}  
\item
$\underline{{\rm PSup}}$ denotes the category of complete
lattices $L$ with morphisms $g:{\cal P}(L_1)\to{\cal P}(L_2)$
that preserve unions and satisfy both eq.(\ref{eq:comm}) and 
\par\medskip\noindent 
\centerline{$
\bigvee_1(A)=\bigvee_1(B)\,\Longrightarrow\,\bigvee_2(g(A))=
\bigvee_2(g(B))\,;$}
\par\medskip\noindent
\item
$\underline{{\rm DCHeyt}}$ denotes the category of complete
Heyting algebras $H$ equipped with a disjunctive join dense
closure
${\cal R}:H\to H$ with morphisms ${h:H_1\to H_2}$ that
preserve suprema and satisfy 
$$ {\cal R}_1(A)={\cal R}_1(B)\,\Longrightarrow\,{\cal
R}_2(h(A))={\cal R}_2(h(B))\,;$$
\item $F:\underline{{\rm PSup}}\to\underline{{\rm
DCHeyt}}:L\mapsto\bigl({\cal DI}(L),{\cal R}_{\cal
DI(L)}\bigr)\,,\,g(-)\mapsto{\cal C}_2\bigl(g(-)\bigr)$\,;
\item $G:\underline{{\rm DCHeyt}}\to\underline{{\rm
Sup}}:(H,{\cal R})\mapsto{\cal R}(H)\,,\,h(-)\mapsto{\cal
R}_2\bigl(h(-)\bigr)$\,;
\item $H:\underline{{\rm PSup}}\to\underline{{\rm
Sup}}:L\mapsto L\,,\,g(-)\mapsto
\bigvee_2\Bigl(g\bigl(\downarrow\!(-)\bigr)\Bigr)$\,.
\end{itemize} Note that the object correspondences of $F$ and 
$G$ are indeed those of Coecke (2002a \S 3)\,, Definition
1\,, and in particular that $G:\underline{{\rm
DCHeyt}}\to\underline{{\rm Sup}}$ is a full quantaloidal
morphism but not an equivalence. This fact will constitute
the core of the discussion below.

\bigskip\noindent
{\bf 4. DISCUSSION}      

\medskip\noindent
Recalling that we mentioned in Coecke (2002a) that
``via physical and logical considerations we
rediscover a purely mathematical result by Bruns and Lakser (1970)
on injective hulls of meet-semilattices'', it is then via these
considerations on state transitions and property
transitions as the morphisms equipping operational resolutions that
the different underlying motivations in Bruns and Lakser
(1970) and Coecke (2002a) reveal themselves explicitly in a formal
way.  Indeed, in
Coecke (2002a) we skipped any consideration on morphisms by only
considering object equivalences,
which only requires specification of isomorphisms.  Obviously, there
are different canonical ways
to extend such an object equivalence categorically, e.g., via a
pointwise lift of the {\it chosen}
morphisms of the complete lattices to the
corresponding complete Heyting algebras of distributive
ideals,\footnote{Obviously provided that the chosen morphisms of the
complete lattices are
$\underline{{\rm Set}}$-concrete.} where as well {\it meet}-morphisms, {\it
inf}-morphisms and {\it
sup}-morphisms are candidate morphisms for the complete lattices.
 From this perspective, in Bruns and Lakser (1970) the choice of
morphisms, i.e., {\it meet}-morphisms, is such that it establishes distributive hulls
as injective hulls. It was moreover noted in Harding (1999) that when
defining {\sl distributive morphisms} as those  maps that
preserve finite meets, distributive suprema and the underlying sets
that have distributive suprema, then
$\underline{{\rm Frame}}$ (Johnstone 1982) is a full monoreflective
subcategory of the category of meet-semilattices and distributive
morphisms, with as reflector a functor that assigns an object to its
distributive hull. However, in view of the
application of distributive hulls ``towards dynamic quantum logic'',
we canonically obtain a full but not
faithful quantaloidal correspondence between complete lattices and
complete Heyting algebras equipped with an operational
resolution, this as a consequence of our manifestly different choice
of morphisms.  These considerations
obviously indicate many new open problems, as do the open questions
posed in Coecke (2002a) when restated in the
presence of this dynamical setting. We do mention at this point that recently we
constructed a dynamic logic that realizes the semantics presented in this paper as a
true logic with forward and backward implication and corresponding tensors ---
see Coecke (2002b) and Coecke and Smets
(2001) for a general presentation and Smets (2001) for a detailed
development of the atomistic case (i.e., Boolean propositions) --- and, that has the
representation presented in Coecke (2002a) as statical limit\,.  Implementation of
this framework for the particular case of quantum measurements represented by
projectors acting on the underlying Hilbert space can be found in Coecke and Smets
(2001)\,.  In that case we obtain a family of implicative hooks labeled by properties.
It is argued there that the {\sl transition from either classical or 
constructive/intuitionistic logic to quantum logic} entails besides the introduction
of an additional unary connective operational resolution the shift from a
binary connective implication to a ternary connective where two of the arguments have
an ontological  connotation and the third, the new one, an empirical. This second
aspect of the shift from classical or 
constructive/intuitionistic to quantum will then be the one that requires
{\sl orthomodularity} of the underlying lattice of properties as a crucial feature.

\bigskip      
\noindent
{\bf ACKNOWLEDGMENTS}  

\medskip\noindent
I thank David Foulis for proposing the study of propagating actuality sets.
Part of the inspiration for this approach emerged from previous joint
work with David Moore, Sonja Smets and Isar Stubbe.  I thank John Harding, Jan Paseka,
Pedro Resende and Isar Stubbe for additional comments on this paper.

\bigskip\noindent
{\bf REFERENCES}\,\footnote{Preprints and postscript files of
published papers by the
current author can be downloaded at http://www.vub.ac.be/CLEA/Bob/Coecke.html.}

\medskip\noindent
A{\scriptsize BRAMSKY}, S. and V{\scriptsize ICKERS}, S. (1993)
`Quantales, Observational Logic and Process
Semantics', {\it Mathematical Structures in Computer Science} {\bf 3}, 161.
\par
\vspace{2mm}
\par
\noindent
A{\scriptsize MIRA}, H., C{\scriptsize OECKE}, B. and S{\scriptsize
TUBBE}, I. (1998) `How
Quantales Emerge by Introducing Induction  within the Operational
Approach', {\it Helvetica Physica Acta}
{\bf 71}, 554.
\par
\vspace{2mm}
\par
\noindent
B{\scriptsize RUNS}, G. and L{\scriptsize AKSER}, H. (1970)
`Injective Hulls of Semilattices',
{\it Canadian Mathematical Bulletin} {\bf 13}, 115.
\par
\vspace{2mm}
\par
\noindent
B{\scriptsize ORCEUX}, F. (1994) {\it Handbook of Categorical Algebra
I \& II},
Cambridge University Press.
\par
\vspace{2mm}
\par
\noindent
B{\scriptsize ORCEUX}, F. and S{\scriptsize TUBBE}, I. (2000) `Short
Introduction to Enriched Categories', In\,:
B. Coecke, D.J. Moore and A. Wilce, (Eds.), {\it Current Research in
Operational
Quantum Logic: Algebras, Categories and Languages},
pp.167--194, Kluwer Academic Publishers.
\par
\vspace{2mm}
\par
\noindent
C{\scriptsize OECKE}, B. (2000) `Structural
Characterization of Compoundness', {\it International Journal of
Theoretical Physics} {\bf 39}, 581\,;
arXiv: quant-ph/0008054\,.
\par
\vspace{2mm}
\par
\noindent
C{\scriptsize OECKE}, B. (2002a) `Quantum Logic in Intuitionistic  
Perspective', {\it Studia Logica} {\bf 70}, 411; arXiv: math.LO/0011208\,.
\par
\vspace{2mm}
\par
\noindent
C{\scriptsize OECKE}, B. (2002b) `Do we have to Retain Cartesian Closedness in
the Topos-Approaches to Quantum Theory, and, Quantum Gravity ?, Preprint.
\par
\vspace{2mm}
\par
\noindent
C{\scriptsize OECKE}, B. and M{\scriptsize OORE}, D.J. (2000)
`Operational Galois Adjunctions', In\,:
B. Coecke, D.J. Moore and A. Wilce, (Eds.), {\it Current Research in
Operational
Quantum Logic: Algebras, Categories and Languages},
pp.195--218, Kluwer Academic Publishers; arXiv:
quant-ph/0008021.
\par
\vspace{2mm}  
\par
\noindent
C{\scriptsize OECKE}, B., M{\scriptsize OORE}, D.J. and S{\scriptsize
TUBBE}, I. (2001)
`Quantaloids Describing Causation and Propagation for Physical
Properties', {\it Foundations of Physics Letters} {\bf 14}, 133; arXiv:
quant-ph/0009100.
\par
\vspace{2mm}
\par
\noindent
C{\scriptsize OECKE}, B. and S{\scriptsize METS}, S. (2000) `A 
Logical Description
for Perfect Measurements', {\it International Journal of Theoretical
Physics} {\bf 39}, 591\,; arXiv:
quant-ph/0008017.
\par
\vspace{2mm}
\par
\noindent
C{\scriptsize OECKE}, B.  and S{\scriptsize METS}, S. (2001) `The Sasaki-Hook is not
a [Static] Implicative Connective but Induces a Backward [in Time] Dynamic One that
Assigns Causes', Paper submitted to {\it International Journal of
Theoretical Physics\,} for the proceedings of IQSA V, Cesena, Italy, April 2001;
arXiv:quant-ph/0111076\,.
\par
\vspace{2mm}
\par
\noindent
C{\scriptsize OECKE}, B. and S{\scriptsize TUBBE}, I. (1999)
`Operational Resolutions and State
Transitions in a Categorical Setting', {\it Foundations of Physics
Letters} {\bf 12}, 29\,; arXiv:
quant-ph/0008020.
\par
\vspace{2mm}
\par
\noindent
D{\scriptsize ANIEL}, W. (1989)
`Axiomatic Descrition of Irreversable and Reversable Evolution of a
Physical System', {\it Helvetica Physica Acta} {\bf
62}, 941.
\par
\vspace{2mm}
\par
\noindent
F{\scriptsize AURE}, Cl.-A. and F{\scriptsize R\"OLICHER}, A. (1993)
`Morphisms of Projective Geometries and of Corresponding Lattices',
{\em Geometriae Dedicata} {\bf 47}, 25.
\par
\vspace{2mm}
\par
\noindent
F{\scriptsize AURE}, Cl.-A. and F{\scriptsize R\"OLICHER}, A. (1994)
`Morphisms of Projective Geometries and Semilinear Maps', {\em
Geometriae Dedicata} {\bf 53}, 237.
\par
\vspace{2mm}
\par
\noindent
F{\scriptsize AURE}, Cl.-A., M{\scriptsize OORE}, D.J. and
P{\scriptsize IRON}, C. (1995)
`Deterministic Evolutions and Schr\"odinger Flows', {\it Helvetica
Physica Acta} {\bf 68}, 150.
\par
\vspace{2mm}
\par
\noindent
H{\scriptsize ARDING}, J. (1999) Private communication.
\par
\vspace{2mm}
\par
\noindent
J{\scriptsize OHNSTONE}, P.T. (1982) {\it Stone Spaces}, Cambridge
University Press.
\par
\vspace{2mm}
\par
\noindent
K{\scriptsize ALMBACH}, G. (1983) {\it Orthomodular Lattices}, Academic Press.
\par
\vspace{2mm}
\par
\noindent
P{\scriptsize IRON}, C. (1976) {\it Foundations of Quantum Physics},
W.A. Benjamin, Inc.
\par
\vspace{2mm}
\par
\noindent
P{\scriptsize OOL}, J.C.T. (1968) `Baer $^*$-Semigroups and the Logic
of Quantum Mechanics', {\it Communications
in Mathematical Physics} {\bf 9}, 118.
\par
\vspace{2mm}
\par
\noindent
R{\scriptsize ESENDE}, P. (2000) `Quantales and Observational Semantics', In\,:
B. Coecke, D.J. Moore and A. Wilce, (Eds.), {\it Current Research in
Operational
Quantum Logic: Algebras, Categories and Languages},
pp.263--288, Kluwer Academic Publishers.
\par
\vspace{2mm} 
\par
\noindent
R{\scriptsize OSENTHAL}, K.I. (1991) `Free Quantaloids', {\it Journal
of Pure and Applied Algebra} {\bf 77}, 67.
\par
\vspace{2mm}
\par
\noindent
S{\scriptsize METS}, S. (2001): `The Logic of Physical Properties in Static and
Dynamic Perspective', PhD-thesis, Free University of Brussels.
\par
\vspace{2mm} 
\par
\noindent
S{\scriptsize OURBRON}, S. (2000) `A Note on Causal Duality', {\it
Foundations of Physics Letters} {\bf 13}, 357.

\end{document}